\documentclass[letterpaper, 10 pt, conference]{ieeeconf}  

\IEEEoverridecommandlockouts                              

\usepackage{amsfonts,bm}
\usepackage{amssymb,amsmath}
\usepackage{graphicx}
\usepackage{cite}
\usepackage{hyperref}
\usepackage{algorithm2e}
\usepackage{algpseudocode}
\RestyleAlgo{ruled}
\usepackage{ifthen}     
\usepackage{xcolor}

\renewcommand{\epsilon}{\varepsilon}
\renewcommand{\phi}{\varphi}

\renewcommand{\d}{\,\mathrm{d}}
\newcommand{\id}{\mathop{\bf id}}

\newtheorem{proposition}{Proposition}

\newtheorem{example}{Example}

\renewcommand{\epsilon}{\varepsilon}
\renewcommand{\phi}{\varphi}

\renewcommand{\d}{\,d}

\newcommand{\R}{\mathbb{R}}

\title{\LARGE \bf
From Few-Shot Optimal Control 
to Few-Shot 
Learning$^*$
}

\author{Roman Chertovskih, Nikolay Pogodaev, Maxim Staritsyn, and A. Pedro Aguiar, \IEEEmembership{Senior Member, IEEE}
\thanks{The authors acknowledge the financial support of the Foundation for Science and Technology (FCT, Portugal) in the framework of the Associated Laboratory ARISE (LA/P/0112/2020), R\&D Unit SYSTEC (base UIDB/00147/2020 and programmatic UIDP/00147/2020 funds) and the bilateral project 2023.09597.CBM. A part of the computations was carried out under FCT computational projects 2023.10674.CPCA and 2024.07885.CPCA.
}
\thanks{Roman Chertovskih, Maxim Staritsyn, and A. Pedro Aguiar are with Research Center for Systems and Technologies (SYSTEC-ARISE), Faculdade de Engenharia, Universidade do Porto, Rua Dr. Roberto Frias, s/n 4200-465, Porto, Portugal (e-mails: roman@fe.up.pt, staritsyn@fe.up.pt, apra@fe.up.pt). Nikolay Pogodaev is with Dipartimento di Matematica ``Tullio Levi-Civita'' (DM), University of Padova, Via Trieste, 63 - 35121 Padova, Italy (e-mail: nickpogo@gmail.com)
}%
}

\begin{document}
\maketitle
\thispagestyle{empty}
\pagestyle{empty}
\begin{abstract}
We present an approach to solving unconstrained nonlinear optimal control problems for a broad class of dynamical systems. This approach involves lifting the nonlinear problem to a linear ``super-problem'' on a dual Banach space, followed by a non-standard ``exact'' variational analysis,~--- culminating in a descent method that achieves rapid convergence with minimal iterations. We investigate the applicability of this framework to mean-field control and discuss its perspectives for the analysis of information propagation in self-interacting neural networks.
\end{abstract}

\section{INTRODUCTION}

This paper develops and numerically validates a monotone descent method for the ``few-shot'' solution of unconstrained nonlinear optimal control problems. The method arises from a non-standard variational analysis of state-linear problems in abstract spaces and extends to a broad class of nonlinear problems via ``linearization by duality''. Applications include optimal mean field control (MFC)~--- an active research area in applied mathematics with uses in physics, mathematical biology, and socio-economic systems (see, e.g., \cite{cristiani2014multiscale,carrillo2010flocking}).

In the last few years, growing interest in this field has been sparked by a series of studies establishing conceptual connections between MFC and machine learning (ML) in artificial neural networks (NNs). This line of research was 
originated from the observation that the conventional architecture of residual NNs (ResNets) can be naturally interpreted as a time-discretized ordinary differential equation (ODE)~\cite{Chen2018,Zuazua}, while the training process is naturally formulated as an optimal control problem~\cite{E2019}. This perspective has provided significant control-theoretical and geometrical insights into the mechanisms of ResNets and AI interpretability. 

Notable results include reformulations of the universal approximation property via ensemble controllability~\cite{Cuchiero,Agrachev2021,Zuazua}, shedding new light on some fundamental ML challenges such as overparameterization and overfitting. Additionally, this viewpoint has inspired new training algorithms based on indirect methods from optimal control~\cite{BONNET2023113161,CSP-neur} (we only mention a few works).

More recently, the success of innovative AI architectures, such as Transformers, has further fueled this topic. In particular, the Transformer's self-attention mechanism~\cite{Vaswani2017}, a key component of many state-of-the-art models, has been naturally interpreted as an interaction term in the corresponding mean field dynamics~\cite{Geshkovski2024,bruno2024emergencemetastableclusteringmeanfield,kim2024transformerslearnnonlinearfeatures}.

\subsection{Contribution and Novelty. Organization of the Paper}

This work unifies and significantly generalizes our previous results \cite{CSP-neur,chertovskihOptimalControlDistributed2023,chertovskihOptimalControlDiffusion2024} on indirect numerical methods for optimal control of \emph{linear} continuity and Fokker-Planck equations, guided by the analysis of exact increment formulas ($\infty$-order variational representations) in cost functionals. In the present paper, we demonstrate that our methodology extends naturally to \emph{fully nonlinear} problems by immersing them into certain linear ``super-problems''. 

To establish a foundation, we begin our analysis in Section~\ref{sec:gc} with the classical optimal control problem on a manifold~--- a topic of independent interest. Notably, even in this case, the proposed approach has not been previously introduced in the literature.

The main contribution is the generalization of this framework to the numerical solution of ``control-linear'' MFC problems with broadband (open-loop) control signals, detailed in Section~\ref{sec:mfc}. Additionally, in Section~\ref{sec:mlc}, we explore potential applications of our control-theoretical findings to ML, highlighting several implications for understanding latent dynamics in modern NN architectures.

Adhering to the format of a short communiqu\'{e}, we deliberately sacrifice some mathematical rigor and omit many theoretical details (some of which can be found in~\cite{pogodaev2025superdualitynecessaryoptimalityconditions}). Instead, we focus on computational challenges and numerical illustrations.

\section{THE INTUITION: GEOMETRIC CONTROL}\label{sec:gc}

Consider the unconstrained optimal control problem in Mayer's form, defined over an $n$-dimensional smooth manifold $\mathcal{X} \doteq \mathcal M^n$ on a fixed time horizon $I = [0, T]$:\footnote{All relations involving measurable functions of time are assumed to hold almost everywhere (a.e.) on $t \in I$. For brevity, this assumption will be implicit in the following discussion.}
\begin{align}
  (G) \quad  \min \ & \left\{\ell(x(T)) \colon x(t) = x[u](t), \ u \in \mathcal{U} \doteq L^\infty(I; U)\right\} \nonumber \\
    \text{subject to} \ & \dot{x}(t) = f_t(x(t), u(t)), \ t \in I, \quad x(0) = \mathrm x_0. \label{mODE}
\end{align}
Here, $\mathrm x_0 \in \mathcal X$ and $\ell \colon \mathcal X \to \R$ are given data. $L^\infty(I; U)$ represents the space of essentially bounded measurable maps $I \to \mathbb{R}^m$ with values in a given convex compact subset $U \subset \mathbb{R}^m$.  In line with geometric control theory \cite{agrachevControlTheoryGeometric2004}, we assume the controlled vector field $f$ is expressed as \(f_t \doteq u_i(t) \, f_t^i,\) with given basis functions $f^i \colon I \times \mathcal{X} \to \mathcal{T}\mathcal{X}$, where $\mathcal{T}\mathcal{X}$ denotes the tangent
bundle of $\mathcal{X}$.\footnote{For brevity, we express the dependence on time in functions of several variables by a subscript and adhere to Einstein's summation convention over other repeated indices.} 


\subsection{Challenges}

The dynamic optimization problem \((G)\) is inherently non-convex, even when \( \mathcal{M}^n = \mathbb{R}^n \) and both \( \ell \) and \( f_t^i \) are linear maps, rendering the problem bilinear. This non-convexity introduces significant challenges, particularly for numerical treatment of the problem. Here, indirect methods with potential for global optimization, such as the shooting method (based on Pontryagin’s Maximum Principle, PMP) and dynamic programming, are typically feasible only in low-dimensional settings. Consequently, practitioners often resort to empirical methods without guaranteed results, gradient-based methods with backtracking \cite{borzì2023sequential}, or (semi-) direct algorithms that rely on problem discretization \cite{ross2012review,betts2010practical,kushner}~--- approaches that are susceptible to the curse of dimensionality and generally provide only local solutions.

This complexity intensifies dramatically for problems involving high-dimensional interacting systems and partial differential equations (PDEs). In such cases, the computational burden is primarily dictated by the numerical solution of the control system itself, rather than the control update process. This scenario underscores the importance of minimizing the number of control updates followed by the recomputation in the state or adjoint equations.

\subsection{Summary of the Approach}

In this paper, we propose an indirect, fully deterministic, and monotone approach comprising the following key steps:

\begin{enumerate}
    \item \emph{Lifting:} The nonlinear problem is transformed into a state-linear ``super-problem'' on the dual Banach space by defining a generating family of linear operators.
    \item \emph{Duality argument:} The linearity facilitates a duality-based framework, leading to a natural definition of the adjoint trajectory in the linear super-problem, which serves as a ``super-adjoint'' to the original trajectory.
    \item \emph{Exact variational analysis:} The variation of the functional in the super-problem is explicitly represented in terms of the super-adjoint trajectory, enabling an exact ($\infty$-order) variational analysis of the original problem.
    \item \emph{Descent method:} Utilizing the increment formula, we derive a feedback control mechanism that ensures a monotone decrease in the cost functional, culminating in a descent algorithm with a self-adjusting step size.
\end{enumerate}
In the rest of this paragraph, we provide a detailed exposition of this strategy.

\subsubsection{Lifting}

Fix $u \in \mathcal{U}$ temporarily. The first idea involves reformulating the nonlinear dynamics $t \mapsto x_t = x_t[u]$ in the \emph{finite-dimensional nonlinear} state space $\mathcal{X}$ as a linear evolution $t \mapsto \mu_t \doteq \delta_{x_t}$ over the \emph{infinite-dimensional vector space} $\bm{X}'$, dual to the Banach space $\bm{X} \doteq C_b(\mathcal{X})$ of bounded continuous functions $\mathcal{X}\to \mathbb{R}$. Here, $\delta_{\rm x}$ denotes the Dirac measure concentrated at the point ${\rm x} \in \mathcal{X}$, establishing the continuous embedding $\mathcal{X} \hookrightarrow \bm{X}'$ though the map ${\rm x} \mapsto \delta_{\rm x}$.

This transformation can be viewed as an extension of the mathematical technique introduced in the theory of Koopman operators \cite{koopman1931hamiltonian} and chronological calculus \cite{agrachevControlTheoryGeometric2004}, suggesting the reconstruction of information about $x_t$ from the dynamics $t \mapsto \phi(x_t)$ in a suitable space $\bm{D} \subset \bm{X}$ of test function; this aligns with the standard approach in statistical mechanics, where test functions act as {observables}. 

Take $\bm D = C^1_c(\mathcal X)$. For a.e. \( t \in I \), introduce the linear unbounded operators \( \mathfrak{L}_t: \bm{D} \to \bm{X} \), sharing a common dense domain \( \bm{D} \), as:  
\begin{align}
    \mathfrak{L}_t \phi \doteq \lim_{\epsilon \to 0^+} \frac{\phi\circ \Phi_{t, t+\epsilon} - \phi}{\epsilon}, \quad \phi \in \bm{D},\label{L}
\end{align}
where \((s,t, \mathrm x) \mapsto \Phi_{s,t}(\mathrm x),\) $I \times I \times \mathcal{X} \to \mathcal{X}$, denotes the flow (or propagator \cite{kolokoltsov2019differential}) of the time-dependent vector field \( (t, \mathrm x) \mapsto f_t(\mathrm x) \doteq f_t(\mathrm x, u_t) \), i.e., \( t \mapsto \Phi_{s, t}(\mathrm x) \) represents the solution to the ordinary differential equation (ODE) \eqref{mODE} with the Cauchy condition \( x(s) = \mathrm x \). 

The operators \eqref{L}, known as Lie derivatives of a function (0-tensor) w.r.t. the vector fields \( f_t \) \cite{Lee2013}, serve as infinitesimal generators of the ODE \eqref{mODE} \cite{kolokoltsov2019differential} in the sense that they define the ``derivative of test functions w.r.t. the corresponding dynamical system'':  
\begin{align}
    \frac{\d}{\d t} \phi(x(t)) = (\mathfrak{L}_t \phi)(x(t)).\label{D}
\end{align}
In local coordinates, the action of \( \mathfrak{L}_t \) on a function $\phi$ is represented as
$
\mathfrak{L}_t \phi = \nabla \phi \cdot f_t$, 
where \( \nabla = \nabla_{\rm x} \) stands for the gradient operator w.r.t. \( {\rm x} \in \mathbb{R}^n \).

Denoting by \( \langle \cdot, \cdot \rangle \) the natural duality pairing in \( (\bm{X}', \bm{X}) \), we rewrite \eqref{D} as  
\begin{align}\label{wfl}  
\frac{d}{dt} \langle \mu_t, \phi \rangle = \langle \mu_t, \mathfrak{L}_t \phi \rangle = \langle \mathfrak{L}_t' \mu_t, \phi \rangle, \qquad \forall\, \phi \in \bm{D},  
\end{align}  
where \( \mathfrak{L}_t' \) denotes the formal adjoint of \( \mathfrak{L}_t \), given locally by the divergence operator:
\(
\mathfrak{L}_t' \vartheta \doteq - {\rm div}(f_t \vartheta).
\)
The relation \eqref{wfl} is recognized as the distributional form of the linear PDE  
\begin{align}  
\partial_t \mu_t = \mathfrak{L}_t' \mu_t, \label{mPDE}  
\end{align}  
which represents an abstract form of the classical continuity (Liouville) equation. Under the standard Cauchy–Lipschitz regularity of the vector fields \( f^i \), this equation admits a unique distributional solution \( t \mapsto \mu_t = \mu_t[u] \) for any initial condition \( \mu_0 = \vartheta \) in the spaces \( {\mathcal P}_1(\mathcal{M}^n) \subset \bm X' \) of probability measures on \( \mathcal{M}^n \) with a finite first moment \cite{ambrosioGradientFlowsMetric2005,Pukhlikov2004} (and even in more general cases, see \cite[Remark 2.4]{pogodaev2025superdualitynecessaryoptimalityconditions}). This solution has an explicit \emph{characteristic} representation:  
\[
\mu_t = (\Phi_{0,t})_\sharp \vartheta,
\]  
expressed in terms of the flow \( \Phi \) \cite{ambrosioGradientFlowsMetric2005}. In particular, the unique solution of \eqref{mPDE} with the initial datum \( \mu_0 = \delta_{x_0} \) is given by the map \( t \mapsto \delta_{x(t)} \), which corresponds to a trajectory of the characteristic system \eqref{mODE}.

Finally, by expressing the cost functional as $\ell(x(T)) = \langle \mu_T, \ell \rangle$, the problem $(G)$ is immersed into the following problem in the dual Banach space:
\begin{align}
\min \left\{ \langle \mu_T, \ell \rangle \colon \mu_t = \mu_t[u], \ \mu_0 = \vartheta, \ u \in \mathcal{U} \right\}.\label{LP}
\end{align}
This new problem is equivalent to $(G)$ when $\vartheta = \delta_{x_0}$, while making rigorous sense for a much broader class of ``distributed'' initial data $\vartheta \in \bm{X}'$, as discussed above. 

\subsubsection{Duality argument} A key feature of the resulting ``super-problem'' is its linearity w.r.t. the new state variable $\mu$. This fact enables the standard duality argument, based on the natural formulation of the adjoint equation to \eqref{mPDE}, as detailed below.

Similar to the definition of the operator $\mathfrak{L}_t$, taking a sufficiently regular map $p \colon I \times \mathcal{X} \to \mathbb{R}$, we have, by the chain rule:
\[
\frac{d}{dt} \langle \mu_t, p_t \rangle = \langle \mu_t, \left\{ \partial_t + \mathfrak{L}_t \right\} p \rangle.
\]
Then, specifying $(t, \mathrm x) \mapsto p_t({\rm x})$ as a mild\footnote{In the sense that $p$ can be differentiable w.r.t. $t$ only a.e. over $I$ with a measurable derivative $t \mapsto \partial_t p_t(\mathrm x)$ for a fixed ${\rm x}$.} solution to the backward (non-conservative) transport PDE
\begin{equation}
    \left\{ \partial_t + \mathfrak{L}_t \right\} p = 0, \quad p_T = \ell,\label{adjmPDE}
\end{equation}
we have: $\frac{d}{dt} \langle \mu_t, p_t \rangle = 0$, implying that
\[
    \langle \mu_t, p_t \rangle \equiv \langle \mu_T, p_T \rangle = \langle \mu_T, \ell \rangle \qquad \forall t \in I.
\]

Remark that a solution to \eqref{adjmPDE} also admits an explicit \emph{characteristic} representation \cite{ambrosioGradientFlowsMetric2005}:
\begin{equation}
p_t = \ell \circ \Phi_{t,T}, \quad t \in I.\label{prep}
\end{equation}

\subsubsection{Exact variational analysis}

Return to the dependence on the control function, inherited by the flow $\Phi = \Phi[u]$ and the generators $\mathfrak{L}_t = \mathfrak{L}_t[u]$, and recall that, due to the control-linear structure of the driving field, we can decompose the operator $\mathfrak{L}_t$ as:
$
\mathfrak{L}_t[u] = u_i(t) \, \mathfrak{A}_t^i,
\)
where operators $\mathfrak{A}_t^i$ generate the flow $\Phi^i$ of the corresponding direction field $f^i$. The key ingredient of our analysis is the following \emph{exact} cost-increment formula in the nonlinear problem $(G)$.
\begin{proposition}
Let $\bar u, u \in \mathcal U$ be arbitrary controls with the corresponding primal and dual states $\mu = \mu[u]$, $\bar \mu = \mu[\bar u]$ and $\bar p = p[\bar u]$. 
Then, the following representation holds:
\begin{align}
\ell(x(T))  - \ell(\bar x(T)) = \int_I \left(u_i(t) - \bar u_i(t)\right) \,\mathfrak A_t^i\bar p_t\big|_{x(t)} \d t.\label{if1}
\end{align}
\end{proposition}
\smallskip
\noindent \begin{proof}
Using the duality argument: \(\langle \bar \mu_t, \bar p_t \rangle \equiv \langle \bar \mu_T, \ell \rangle\), and the initial conditions: $\mu_0 =\bar \mu_0 = \vartheta$, 
we represent
\begin{align*}
     &\langle \mu_T - \bar{\mu}_T, \bar p_T\rangle - \underbrace{\langle\mu_0 - \bar{\mu}_0, \bar p_0\rangle}_\text{$ = 0$}\nonumber= \langle \mu_T, \bar{p}_T\rangle - \langle \mu_0, \bar p_0\rangle\\
& = \int_I \frac{\d}{\d{t}} \langle \mu_t, \bar{p}_t\rangle \d t= \int_I \left(u_i(t) - \bar u_i(t)\right)\langle \mu_t, \mathfrak A_t^i\bar p_t\rangle {\d} t.
\end{align*}%
Setting in this expression $\mu_t = \delta_{x(t)}$, and recalling that $\bar p_T = \ell$, we come to the desired expression \eqref{if1}.
\end{proof}

We now demonstrate how this formula can be used to test optimality. Let \( \bar{u} \) serve as a reference (given or computed) control, while the control \( u \) is sought to provide descent: \( \ell(x(T)) < \ell(\bar{x}(T)) \). A straightforward approach is to choose \( u \) as a pointwise minimizer in the integrand:  
\begin{equation}
    \label{u-max}
    u_i(t) \left(\mathfrak{A}_t^i \bar{p}_t\right)(x(t)) = \min_{(\upsilon_i) \in U} \upsilon_i \left(\mathfrak{A}_t^i \bar{p}_t\right)(x(t)),
\end{equation}
which automatically implies 
\begin{align}
\ell(x(T)) \leq \ell(\bar{x}(T)).\label{nas} 
\end{align}
However, it must be emphasized that this expression involves a feedback loop, as \( x = x[u] \) depends on the yet-to-be-determined input \( u \) in an \emph{operatorial} manner. 
Below, we will show that this equation can be numerically resolved by 
a standard sample‐and‐hold approach. 

\subsubsection{Descent method}
We now turn to the numerical (approximate) solution of problem~$(G)$. 

For simplicity, let $U = [-1,1]^m$, in which case the $i$th control component is exactly responsible for ``switching on/off'' the corresponding direction field~$f^i$.

The corresponding solution to the problem~\eqref{u-max}, which turns out to be separable w.r.t.~$i$, is given by
\[
u_i(t) = - {\rm sign}\left(\left(\mathfrak A_t^i \bar p_t\right)(x(t))\right),
\]
where we set ${\rm sign}(0) = 0$ for numerical stability.\footnote{Standard approaches~\cite{borzì2023sequential}, also followed by us in \cite{chertovskihOptimalControlDiffusion2024,chertovskihOptimalControlDistributed2023}, introduce a small penalization in the $L^2$-norm of~$u$. However, this step significantly affects the structure of solutions. Here, we aim to obtain piecewise constant controls with a small number of switches, preserving the original problem structure and satisfying practical requirements in control design.}

Recalling the representation of~$\mathfrak A_t^i \bar p_t$ via the flows~$\Phi^i$, we approximate:
\begin{align*}
    \mathfrak A_t^i \bar p_t &\approx \frac{\bar p_t(\Phi^i_{t, t+\epsilon}) - \bar p_t}{\epsilon}
    = \frac{\ell\circ\bar \Phi_{t, T}\circ \Phi^i_{t, t+\epsilon} - \ell\circ\bar \Phi_{t, T}}{\epsilon}.
\end{align*}

Since the factor~$1/\epsilon$ does not affect the sign of the expression, we obtain an approximate feedback law for synthesizing descent directions from~$\bar u$:
\begin{equation}
    u_{i}^{\epsilon}(t) = {\rm sign}\left\{\ell\circ\bar \Phi_{t, T} - \ell\circ\bar \Phi_{t, T}\circ \Phi^i_{t, t+\epsilon}\right\}\Big|_{x(t)}.\label{u}
\end{equation}
Here, the expression in brackets represents the gain obtained by a short-time activation of the direction field~$\pm f^i$. 

This representation suggests a simple mechanism for synthesizing descent controls in the Krasovskii-Subbotin fashion \cite{krasovskii2011game}. Introducing a short-time planning horizon~$h$, we define:
\begin{align}
\mathrm u_i\big|_{[0,h]} = {\rm sign}\left\{\ell(\bar \Phi_{0,T}(\mathrm x_0)) - \ell(\bar \Phi_{0,T}(\Phi^i_{0,\epsilon}(\mathrm x_0)))\right\},\label{eu}
\end{align}
which involves solving the ODE~\eqref{mODE} $(m+2)$ times:
\begin{enumerate}
\item First, integrate the system on the entire interval~$I$ with initial condition~${\rm x}_0\) and reference control~$\bar u$ to compute the reference terminal point~$\bar{\rm x}_T$ and the corresponding nominal cost~$\ell(\bar{\rm x}_T)$.
\item Second, solve~\eqref{mODE} on the short-time interval~$[0,\epsilon]$ with the same initial data and control~$(u \equiv e^i)$, obtaining the point~$\mathrm x^i_{h,\epsilon}=\Phi^i_{0,\epsilon}(\mathrm x_0)$; here, $e^i$~denotes the $i$th basis unit vector in~$\mathbb R^m$.
\item Solve~\eqref{mODE} again on~$I$ with reference control~$\bar u$ and initial data~$\mathrm x^i_{h,\epsilon}$ to compute the corrected cost~$\ell(\mathrm x^i_{h,\epsilon})$.
\end{enumerate}

By setting~$h=T/N$ and applying this strategy iteratively on the moving window~$[t,t+h]$, we obtain Algorithm~\ref{algo}. 
This method consists of repeatedly recomputing a constant control over a sliding short-time horizon. It can be shown~\cite{krasovskii2011game} that, as~$\max\{\epsilon, 1/N\}\to 0$, the resulting piecewise constant controls converge\footnote{up to a subsequence in the weak* topology~$\sigma(L^\infty,L^1)$ of the space~$\mathcal U$} to a measurable solution of the \emph{operator equation} \eqref{u-max} ensuring the non-ascendancy property~\eqref{nas}. 

Although each iteration of this method requires solving the state equation $(m+2)N$ times, in many applications exhibiting \emph{sufficiently strong controllability}, the number of sampling points~$N$ and outer iterations~$N_{\mathrm{iter}}$, needed to achieve an acceptable solution, is typically quite small. For instance, in the examples below, we have~$N=3$ and~$N_{\mathrm{iter}}=1$. This results in a significant computational advantage compared to gradient-type algorithms.  Additionally, the optimized control is constructed as a \emph{piecewise-constant} function with \emph{a~priori fixed} instants of potential switches, a desirable feature that aligns with the practice of control engineering.

\begin{algorithm}\label{algo}

\small 
\SetAlgoLined
\KwIn{Initial state \(\mathrm x_0\), final time \(T\), \# MPC intervals \(N\), perturbation \(\varepsilon\), \# iterations \(N_{\mathrm{iter}}\), initial control \(\bar u\)}
\KwOut{Control policy \(u\) and state trajectory \(x\)}
\For{\(k=1,\dots,N_{\mathrm{iter}}\)}{

  \tcp{Initialize synthesis} \(t\gets0\), \(\mathrm x\gets \mathrm x_0\)\; 
 
  \While{\(t<T\)}{

    \tcp{Compute reference cost from $\rm x$}
    \(\bar{\rm y}^t\gets \Phi_{t,T}[u = \bar u](\mathrm x)\), \quad \(\bar J^t=\ell(\bar{\rm y}^t)\)\;
    
    \For{\(i=1,\dots,m\)}{
    
      \tcp{Perturb \(i\)th control direction}
      \({\rm y}^{t,\epsilon}_{i}\gets \Phi_{t,t+\epsilon}[u \equiv e^i](\mathrm x)\)\; 

       \tcp{Compute perturbed cost from $\rm x$}
      
      \(\mathrm z^{t, \varepsilon}_{i}\gets \Phi_{t,T}[u = \bar u]({\mathrm y}^{t,\epsilon}_{i})\), \quad \(J_{i, \epsilon}^t=\ell(\mathrm z_i^{t,\epsilon})\)\;
    }
    \tcp{Update control and propagate}
 
    \(\mathrm u_i^t \gets {\rm sign}(\bar J^t-J_{i, \epsilon}^t)\)\;
    
    \(\mathrm x\gets \Phi_{t,t+h}[u \equiv (\mathrm u_i^t)](\mathrm x)\), \(t\gets \min(t+1/N,\,T)\)\;

  }
  \tcp{Update baseline control }
  \(\bar u\gets u\)\;
}
\Return Final control policy and trajectory
\caption{(metacode)}
\end{algorithm}

\subsection{Generalizations} 

A brief inspection of the previous section reveals that the presented analysis neither relied on the local (manifold) structure of the state space~\(\mathcal{X}\), nor explicitly used the detailed knowledge of the vector field. Instead, the essential ingredients were the definitions and regularity properties of the flow~\(\Phi\) and the generating family~\((\bm D, \mathfrak L_t)\).

This observation naturally suggests extending the approach to a broader class of control systems~\((\mathcal{X}, \Phi, \mathcal{U})\), including those in which~\(\mathcal{X}\) is infinite-dimensional, and~\(\Phi = \Phi[u]\) denotes the solution operator associated with a nonlinear distributed system.  
We illustrate this extension below by redefining the objects~\((\mathcal X, \Phi)\) and~\((\mathfrak L_t, \bm D)\) within the framework of mean-field control~(MFC).

\section{EXTENSION: MEAN-FIELD CONTROL}\label{sec:mfc}

Let now $\mathcal X = \mathcal P_2 \doteq (\mathcal P_2(\mathcal M), W_2)$ be the metric space of probability measures with a finite second moment over a manifold $\mathcal M = \mathcal M^n$, equipped with the 2-Kantorovich (Wasserstein) distance~\cite{ambrosioGradientFlowsMetric2005} and Otto tangent bundle~\cite{Ott2001}. 

Consider~$t \mapsto \bm \Phi_{s,t}[u](\vartheta)$ as the distributional solution to the controlled nonlocal continuity equation:
\begin{equation}
    \partial_t \mu_{t} + \nabla \cdot (\bm f_t(\mu_t,u_t)\,\mu_{t}) = 0,\quad \mu_s = \vartheta,\label{PDE}
\end{equation} 
where \(\bm f_t \doteq u_i(t)\,\bm f_t^i\), and~\(\bm f^i\colon I\times \mathcal M\times \mathcal P_2\to \mathcal T\mathcal M\) are \emph{nonlocal} driving fields satisfying standard regularity hypotheses~\cite{pogodaev2025superdualitynecessaryoptimalityconditions}, and $u$ belongs to the same class~$\mathcal U$.\footnote{This implies that the control~$u_t$ acts uniformly on the entire distribution~$\mu_t$, placing our model within the context of ``ensemble control''~\cite{Brockett2007}.}

Fixing $\vartheta\in\mathcal P_2$, we redefine~$\mu_t[u]\doteq \Phi_{0,t}[u](\vartheta)$ and introduce the fully nonlinear version of the problem~\eqref{LP}:
\begin{align}
 (MF)\quad \min &\left\{\bm \ell(\mu_T)\colon \mu_t = \mu_t[u],\ u\in\mathcal U\right\}.
\end{align}

Analogously to the classical case, this problem can be lifted into a linear super-problem at a higher level of abstraction. Roughly speaking, the nonlinear equation~\eqref{PDE} now plays the role of a ``characteristic equation'' associated with an abstract linear ``PDE'' defined on the space of ``measures over probability measures.''

To formalize this step rigorously, one needs to accurately identify the generating family~$(\bm D,\mathfrak L_t\sim\{\mathfrak A_t^i\})$ associated with the PDE~\eqref{PDE}, including a suitable class of test functionals\footnote{We emphasize the distinction between test functionals, denoted by bold letters~$\bm \phi$, and the ``lower-level'' test functions~$\phi\in C_c^1(\mathcal M^n)$ used in defining distributional solutions of PDE~\eqref{PDE}.}~$\bm \phi\colon\mathcal X\to\mathbb R$ and the generators of flows~$\bm\Phi^i$ corresponding to the basis nonlocal fields~$\bm f^i$.

A natural choice for~$\bm D$ is the set of \emph{intrinsically} differentiable \cite{CardMaster2019,cardaliaguetAnalysisSpaceMeasures2019} functionals of probability measures with~$\mathcal C^{1,1}$ regularity~\cite{chertovskihOptimalControlNonlocal2023}. For~$\bm \phi\in\bm D$, formal calculation of~the pathwise derivative $\frac{d}{dt}\bm \phi\circ\bm\Phi_{0,t}$, performed according to~\cite{CardMaster2019}, yields the explicit representation:
\[
   \mathfrak A_t^i \bm\phi\doteq \big\langle \bm\nabla_\vartheta\bm\phi,\bm f_t^i(\cdot,\vartheta)\big\rangle_\vartheta,
\]
where~$\bm\nabla_\vartheta$ denotes the intrinsic gradient, and~$\big\langle\cdot,\cdot\big\rangle_\vartheta$ is the scalar product in the weighted space~$L^2_\vartheta$.\footnote{Due to space limitations, we omit technical details and refer the reader to~\cite{cardaliaguetAnalysisSpaceMeasures2019} for a comprehensive introduction to analysis and geometry on~$\mathcal P_2$, and to~\cite{Villani2009,ambrosioGradientFlowsMetric2005} for thorough treatments.}

With these definitions in place, we reproduce the increment formula~\eqref{if1} using the family of functionals
\[
\bar{\bm p}_t = \bm\ell\circ\bm\Phi_{t,T}[\bar u],\quad t\in I,
\]
satisfying an abstract transport equation analogous to~\eqref{adjmPDE} within the updated space of functionals~$\bm D$. This yields a computational method directly analogous to Algorithm~\ref{algo}.

We now illustrate the efficiency of this method with a simple yet insightful benchmark problem presented in~\cite{chertovskihOptimalControlNonlocal2023}.

\begin{example}[Kuramoto model]\label{KuramE}
Consider the setting with $\mathcal{M} = \mathbb{S}^1$, $m=2$, the objective functional~$\bm{\ell}(\mu) = \langle \mu, F \rangle$, where~$F({\mathrm x}) = 1 - \cos(\rm x - \pi)$, and the basis nonlocal vector fields defined as:
\begin{equation*}
    f_1 \equiv 1,\quad f_2(\mathrm x,\mu) = \int_{0}^{2\pi}\sin(\mathrm y - \mathrm x)\,\mathrm{d}\mu(\mathrm y).
\end{equation*}
In other words, the task is to concentrate the terminal density around a given point~$\check{\mathrm x}=\pi$, which models a phase synchronization problem in the Kuramoto-type model~\cite{Kuramoto1975}.

In~\cite{chertovskihOptimalControlNonlocal2023}, this problem was addressed via an indirect numerical algorithm with backtracking, based on Pontryagin’s maximum principle. For the absolutely continuous initial distribution~$\vartheta(\mathrm{d} \mathrm x)=\rho_0(\mathrm x)\,\mathrm{d} \mathrm x$ with
\[
\rho_0(\mathrm x)=\frac{2+\sin \mathrm x +0.8\cos 2 \mathrm x-0.2\sin 2 \mathrm x}{4\pi},
\]
an acceptable solution was obtained after 29 iterations, involving the solution of the primal and adjoint PDEs \emph{188 and 29 times}, respectively.

For the same initial data, identical PDE discretization parameters, and perturbation factor~$\epsilon=0.1$, Algorithm~\ref{algo} achieved a qualitatively similar result after just a single iteration, involving four control switches (see the discussion below), at the cost of solving the primal PDE \emph{12~times}.

We also considered the problem of optimally matching a given target density profile
\(
\hat{\rho}(x)=\exp(-0.5(x-\pi)^2),
\)
by setting~$F(x)=|\rho_T(x)-\hat{\rho}(x)|^2$. This problem was also solved in a single iteration, reducing the cost from~$0.14$ to~$0.007$ with only three control switches (see Fig.~\ref{kuram}). 
Notably, the method presented in~\cite{chertovskihOptimalControlNonlocal2023} is \emph{not} directly \emph{applicable} to this scenario since it relies on computing the intrinsic derivative of the cost functional, which is challenging when the cost explicitly depends on the density function.

\begin{figure}
    \centering

\vspace{0.15cm}
    
    \includegraphics[width=0.3\textwidth]{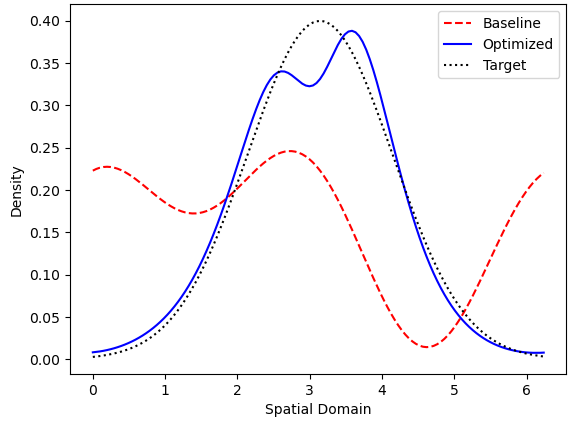}
    \caption{Example~\ref{KuramE}: Baseline, $\rho_0$, and optimized, $\rho_T$, densities together with its target profile $\hat{\rho}$ on the flat  representation  of $\mathbb S^1$, $\mathrm{x} \in [0, 2\pi]$.}
    \label{kuram}
\end{figure}

\end{example}

\section{ML PERSPECTIVES}\label{sec:mlc}

In the last few years, MFC has emerged as a powerful framework for analyzing self-organization phenomena in large-scale NNs, particularly Transformers 
 \cite{Geshkovski2024,kim2024transformerslearnnonlinearfeatures,bruno2024emergencemetastableclusteringmeanfield}. 
 
Transformers process sequences of elementary units, called tokens, which are embedded into a high-dimensional manifold. Inputs are linearly projected through learnable query, key, and value matrices, $(Q,K,V)$, allowing tokens to dynamically attend to each other via the so-called \textit{self-attention} mechanism~\cite{Vaswani2017}. 

As recognized in~\cite{Geshkovski2024,bruno2024emergencemetastableclusteringmeanfield}, the iterative token update via self-attention resembles the evolution of interacting particle systems on latent embedding manifolds. Such evolution admits a mean-field PDE representation~\eqref{PDE}, where the nonlocal vector field~$\bm f$ is modulated by the parameters~$(Q,K,V)$. 

A key phenomenon arising from such dynamics is \emph{clustering}, wherein the probability distribution of token embeddings, representing context, converges towards a consensus state. This convergence effectively collapses empirical measures into a delta function, concentrated at a point in the latent space corresponding to the newly predicted token. In practice, this means that all token embeddings become increasingly similar, emphasizing the most relevant information for prediction. This behavior mirrors \emph{synchronization} effects observed in Kuramoto-type models. As such, it provides a challenging benchmark to evaluate our method, offering numerical evidence that attention-like nonlocal interaction fields can effectively aggregate probability distributions of token embeddings around desired points after learning the parameters~$(Q, K, V)$.

Following~\cite{Geshkovski2024,bruno2024emergencemetastableclusteringmeanfield}, we will demonstrate that such aggregation problems are solvable even in a simplified case where $Q_t = K_t \equiv \id$, meaning that the queries and keys are not transformed and attention operates directly on the input embeddings. In this setting, only the value matrix function $t \mapsto V_t$ is optimized as a control $I \to \mathbb R^{2 \times 2}$, preserving the role of self-attention as an adaptive combination mechanism.

\begin{example}[Attention-inspired model]\label{AtteE}
Let the underlying manifold ${\mathcal M}$ be the 2D torus~$\mathbb T^2$. Such a manifold naturally arises in ML tasks involving periodic data (e.g., texture synthesis, generative modeling on cyclic domains) and provides an appealing testbed for attention mechanisms. To reflect the periodic nature of this setting, we define the attention-like field using a von~Mises-type kernel: 
\[
\bm f_t({\rm x}, \vartheta, V_t) = 
\frac{\displaystyle\int \exp\left(\kappa\sum \cos({\rm x}_j - {\rm y}_j)\right) V_t \ {\rm y}\,\mathrm d\vartheta({\rm y})}
{\displaystyle\int\exp\left(\kappa\sum \cos({\rm x}_j - {\rm y}_j)\right)\,\mathrm d\vartheta({\rm y})},
\]
where $\mathrm{x} = (x_1, x_2) \in [0, 2\pi]^2$ are the local coordinates in the flat 
representation of the torus, and $\kappa>0$ is a concentration parameter.

Starting from a given initial distribution, the aim is to aggregate the density around the position $\check{x}=(0,0)$ corresponding to a corner in the fundamental domain of~$\mathbb T^2$. This is achieved by optimizing the control $t \mapsto V_t$ by our approach, using a modest number of switches. In experiments, we set $T=0.5$, $\epsilon=0.1$ and $\kappa = 5$. Satisfactory results were consistently obtained in a single iteration using $N=4$ control intervals, as shown in Fig.~\ref{fig:atte}. The computations can be reproduced using the code available via 
\url{https://github.com/starmaxwell/Horse.git}.
\begin{figure}[h!]

\vspace{0.18cm}

    \centering
  \includegraphics[width=0.46\textwidth]{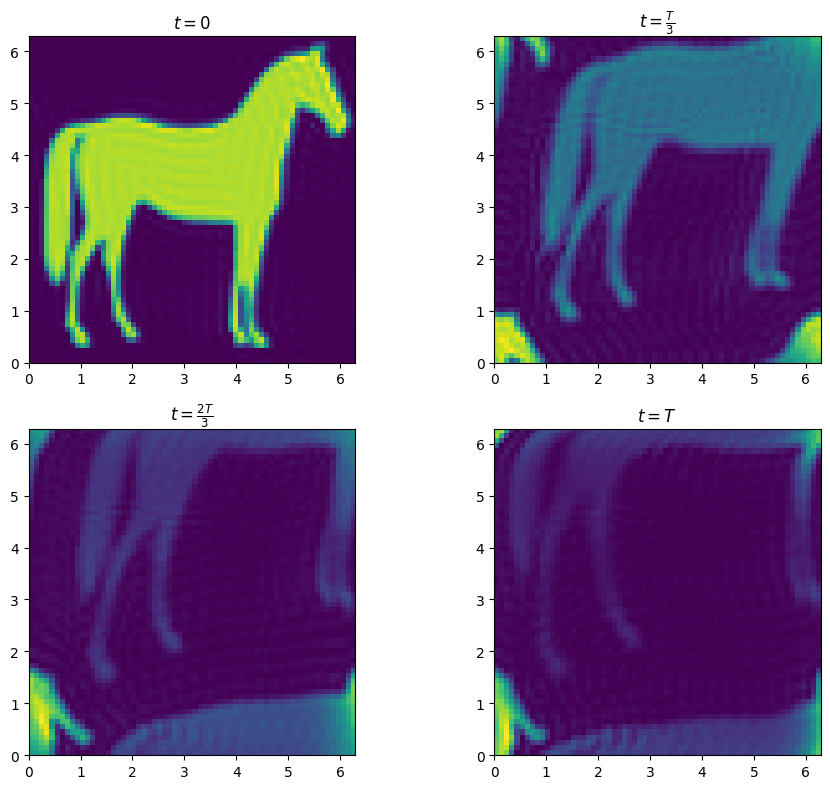}
    \caption{Example~\ref{AtteE}: Optimized density snapshots as heatmaps over the torus' fundamental domain (rendered on a relative color scale).}
    \label{fig:atte}
\end{figure}
\end{example}

\section{CONCLUSION: FEW-SHOT LEARNING}\label{sec:con}

The presented approach is potentially applicable to a broader class of control problems involving dynamical systems on metric and measure spaces, opening several directions for future research.

Applications in ML are also quite promising. Here, the key challenge is the numerical solution of PDEs, typically feasible only in low dimensions, far from those of realistic feature spaces. In view of this limitation, particularly compelling is the application of our optimal control framework to \emph{few-shot learning} (FSL), where classifiers or regressors are trained on limited datasets. Specifically, if each data sample is represented as a probability distribution on a low-dimensional manifold such as $\mathbb{T}^2$, the FSL problem is naturally stated as MFC in a \emph{system} of coupled continuity equations~--- each describing the embedding of a corresponding data sample into the manifold~--- thereby rendering the problem computationally manageable. This perspective aligns with recent developments in PDE-based FSL \cite{Wang2019}. 

\bibliographystyle{IEEEtran}
\bibliography{IEEEabrv,starmax_full}

\end{document}